\documentclass[letterpaper, 10 pt, conference]{ieeeconf}  % Comment this line out if you need a4paper

\IEEEoverridecommandlockouts                              % This command is only needed if 
\usepackage{hyperref}
\usepackage{url}
\usepackage{bbold}
\usepackage{amsfonts}
\usepackage{amsmath,amssymb}
\usepackage{graphicx}
\usepackage{hyperref}
\usepackage{wrapfig}
\usepackage{mathrsfs} 
\usepackage{algorithm,algorithmic}
\usepackage{array}
\usepackage{times}
\usepackage{url}
\usepackage{cite}
\usepackage{upgreek}
\usepackage{float}
\usepackage{color}
\usepackage{bbm}

\newtheorem{theorem}{Theorem}

\newtheorem{definition}{Definition}

\newtheorem{lemma}{Lemma}

\newtheorem{assumption}{Assumption}

\newcommand\numberthis{\addtocounter{equation}{1}\tag{\theequation}}

\newcommand{\R}{\mathbbm{R}}
\newcommand{\e}{\mathbf{e}}
\newcommand{\1}{\mathbbm{1}}
\newcommand{\0}{\mathbf{0}}
\newcommand{\E}{\mathbbm{E}}
\newcommand{\w}{\mathbf{w}}

\newcommand{\y}{\mathbf{y}}

\newcommand{\te}{\hat{\boldsymbol \theta}}

\newcommand{\diag}{\textit{\textbf{diag}}}

\title{\LARGE \bf 
Distributed Estimation of Dynamic Parameters : Regret Analysis
}

\author{Shahin Shahrampour, Alexander Rakhlin and Ali Jadbabaie
\thanks{This work was supported by ONR BRC Program on Decentralized,
Online Optimization.}
\thanks{Shahin Shahrampour is with the Department of Electrical Engineering at Harvard University, Cambridge, MA 02138 USA. (e-mail: shahin@seas.harvard.edu).}
\thanks{Alexander Rakhlin is with the Department of Statistics at the University of Pennsylvania, Philadelphia, PA 19104 USA. (e-mail: rakhlin@wharton.upenn.edu).}
\thanks{Ali Jadbabaie is with the Department of Electrical and Systems Engineering at the University of Pennsylvania, Philadelphia, PA 19104 USA. (email: jadbabai@seas.upenn.edu).}
} %

\begin{document}

\maketitle
\thispagestyle{empty}
\pagestyle{empty}

%%%%%%%%%%%%%%%%%%%%%%%%%%%%%%%%%%%%%%%%%%%%%%%%%%%%%%%%%%%%%%%%%%%%%%%%%%%%%%%%
\begin{abstract}
This paper addresses the estimation of a time-varying parameter in a network. A group of agents sequentially receive noisy signals about the parameter (or moving target), which does not follow any particular dynamics. The parameter is not observable to an individual agent, but it is globally identifiable for the whole network. Viewing the problem with an online optimization lens, we aim to provide the finite-time or non-asymptotic analysis of the problem. To this end, we use a notion of dynamic regret which suits the online, non-stationary nature of the problem. In our setting, dynamic regret can be recognized as a finite-time counterpart of stability in the mean-square sense. We develop a distributed, online algorithm for tracking the moving target. Defining the path-length as the consecutive differences between target locations, we express an upper bound on regret in terms of the path-length of the target and network errors. We further show the consistency of the result with static setting and noiseless observations.
\end{abstract}

\section{Introduction}
Distributed detection, learning and estimation has been a main topic of interest in the past three decades \cite{tsitsiklis1988decentralized,stankovic2011decentralized,drakopoulos2013learning,7349151}. With wide-spread applications in sensor, robotic, economic and social networks, distributed algorithms have received a considerable attention in science and engineering \cite{bullo2009distributed,acemoglu2011bayesian}. In these scenarios, a group of agents aim to learn or estimate the value of a parameter. Each individual agent receives partially informative data about the parameter; however, the global spread of information in the network allows agents to accomplish the task collaboratively. Many of these information aggregation methods use {\it consensus} protocols as a crucial component \cite{jadbabaie2003coordination,olfati2004consensus}. Distributed algorithms are popular for their ability to handle large data sets, low computation burden on agents and robustness to node failures. 

On the other hand, online learning and optimization has been intensively studied in the literature of machine learning \cite{cesa2006prediction,zinkevich2003online}, proving to be a powerful tool to model sequential decision problems. The problem can be modeled as a game between a learner and an adversary. The learner sequentially selects actions, and the adversary reveals the corresponding losses to the learner. The term {\it online} refers to the fact that the learner receives data in a sequential fashion.  The popular performance metric for online algorithms is called {\it regret}. Regret measures the performance of algorithm versus a pre-defined benchmark. For instance, the benchmark could be the best fixed action had the learner known all the losses in advance. In a broad sense, when the benchmark is a fixed sequence, the regret is called {\it static}, whereas a {\it time-varying} benchmark sequence brings forward the notion of {\it dynamic} regret \cite{zinkevich2003online,hall2015online,jadbabaie2015online}.   

The goal of this paper is to develop a {\it distributed, online} algorithm for tracking {\it time-varying} parameters. To this end, a unifying observation is to view the distributed estimation problem as an online optimization. A network of agents aim to track a moving target which is only partially observable to each agent. The parameter can represent a sequence of high-dimensional data where agents only receive a low-dimensional version of it. Therefore, they must communicate with each other to track the target. Each time agents form their estimates of the parameter, the whole network incurs a loss (network loss). We formalize the problem as an instance of distributed, online optimization whose {\it dynamic} regret can be characterized in the general form of
\begin{align*}
\textbf{Regret} = \text{Decentralized Loss - Centralized Loss}.
\end{align*}
In other words, the regret compares the performance of the decentralized algorithm versus its centralized counterpart. 

Our main contribution is to provide the {\it finite-time} or {\it non-asymptotic} analysis of the problem using the notion of regret. In this context, regret can also be interpreted as a finite-time analog of {\it stability} in mean-square sense. In contrast to existing approaches, we do not restrict the dynamics of the parameter, i.e., the target is allowed to move along an arbitrary trajectory. The freedom given to the dynamics broadens the application of the problem. For instance, the parameter can represent the price of a product which does not have a particular pattern, or the location of a vehicle that moves along an unpredictable trajectory. In turn, our estimation update is akin to {\it consensus+innovation} updates in the literature with an important distinction that it does not assume {\it any} dynamics on the parameter. 

We then characterize the behavior of regret in finite-time regime. To this end, we quantify the {\it path-length} in terms of the consecutive differences between target locations, and prove that the regret can be bounded in terms of the path-length of the target. Interestingly, our regret bound interpolates smoothly from static to dynamic setting. That is, when the parameter is static, we can recover the corresponding results with {\it static regret} (with appropriate tuning of problem parameters). We further simplify the bound in two cases: static target and noiseless observation.

\subsection*{Related Work}
There exists a large body of literature on estimation of dynamic parameters. These works mostly assume that the parameter follows a {\it known} dynamics, say, the model can be represented as an LTI system. Perhaps the most classical example is the celebrated Kalman filtering \cite{kalman1960new}. The elegance of Kalman filter has motivated a lot of researchers to investigate the problem in the context of networks. We cannot hope to do the justice to the extensive literature on {\it distributed} Kalman filtering, and refer the avid reader to a series of works on this topic\cite{olfati2005distributed,olfati2007distributed,carli2008distributed,olfati2009kalman}. The works of \cite{acemoglu2008convergence,khan2010connectivity,das2013distributed,shahrampour2013online} are also in the same spirit in which the parameter follows linear dynamics. Acemoglu et al. \cite{acemoglu2008convergence} present a rule-of-thumb learning rule, and provide the {\it asymptotic} behavior of their update. In \cite{khan2010connectivity}, Khan et al. investigate the trade-off between {\it stability} of the linear dynamics and {\it connectivity} of the network. They show that their update can potentially track unstable linear models driven by noise. Similarly, an algorithm consisting of pseudo-innovations is developed in \cite{das2013distributed}. Moreover, the authors of \cite{shahrampour2013online} consider a {\it scalar} linear model, and characterize an explicit expression for the mean-square deviation. Restricting their attention to noiseless case, De et al.  pose an {\it inverse} problem in \cite{de2014learning}. In their model, the parameter can be observed whereas the dynamics is {\it unknown} and must be estimated. On the other hand, Atanasov et al. \cite{atanasov2014joint} propose a distributed, linear estimator in wireless sensor networks which encompasses an auxiliary localization procedure. The authors prove the mean-square {\it consistency} of the joint localization and target estimation algorithm. Moreover, Sayin et al. \cite{sayin2015optimal} present a distributed, online algorithm for state estimation applicable to Big Data. The authors of \cite{simonetto2015decentralized} develop an optimization-based, prediction-correction method for tracking the trajectory of a moving target. Finally, our work is also related to \cite{van2016learning}, where the problem is considered in the context of social networks.

\subsection*{Organization}

The rest of this paper is organized as follows. In the next section we introduce the notation, and formalize the estimation problem using the notion of regret. In Section \ref{Main Results}, we present our technical results and their consequences. Notably, we demonstrate a path-length bound on regret. Section \ref{Concluding Remarks} consists of concluding remarks, and in Section \ref{Appendix : Proofs} we provide the proof of our technical results. 

\section{Problem Formulation}\label{Problem Formulation}

\subsection{Nomenclature} 
We adhere to the following notation in the exposition of our results:
\begin{center}
  \begin{tabular}{| c || l | }
    \hline
     $[n]$ &  The set $\{1,2,...,n\}$ for any integer $n$ \\ \hline
     $x^\top$ & Transpose of the vector $x$ \\ \hline
     $I_d$ &  Identity matrix of size $d \times d$ \\ \hline
     $\|\cdot\|$ &  The spectral norm operator \\ \hline
     $\1_d$ &  Vector of all ones with dimension $d$ \\ \hline
     $\0_d$ &  Vector of all zeros with dimension $d$ \\ \hline
     $\lambda_i(P)$ &  The $i$-th largest eigenvalue of matrix $P$ \\ \hline
     $\rho(P)$ & The spectral radius of matrix $P$ \\ \hline
  \end{tabular}
\end{center}
Throughout, the variables in bold are random, and all the vectors are in column form. When referring to vector of all ones or zeros, sometimes the dimension is omitted and can be inferred from the context.

\subsection{Observation Model and Regret Definition}
Our objective is to track a {\it time-varying} parameter that is not confined to a certain dynamics. The parameter could represent the location of a moving target following an arbitrary trajectory. We denote the parameter by $\theta_t \in \R^{d}$ at each time $t\in [T]$. A network of $n$ agents collaborate with each other to track the parameter which is only {\it partially} observable to each agent. More formally, the observation model for agent $i \in [n]$ can be expressed in the following form
\begin{align}\label{model}
\y_{i,t}=H_i \theta_t + \w_{i,t},
\end{align}
where $H_i \in  \R^{m_i \times d}$ and $\w_{i,t} \in \R^{m_i}$ denote the observation matrix and noise, respectively. For instance, $\theta_t$ can be a high-dimensional data where agents only observe a low-dimensional version of that, i.e. $m_i$ might be much smaller than $d$. We respect the standard assumption of zero-mean and finite-covariance for the observation noise. That is, we assume that $\E\left[\w_{i,t}\right]=\0$ and $\E\left[\w_{i,t}^\top\w_{i,t}\right]=W_i$, respectively. We further assume the independence over time and space; that is, $\E\left[\w_{i,t}^\top\w_{j,s}\right]=0$ for $i\neq j$ or $s\neq t$.

Note that when the parameter is {\it static}, i.e. when $\theta_t=\theta$, we recover the classical distributed estimation problem (see e.g. \cite{stankovic2011decentralized}). Though not locally observable to each agent, we assume that the parameter is {\it globally} identifiable from the standpoint of all agents in the network.
\begin{assumption}\label{A1}
The sequence $\{\theta_t\}_{t=1}^T$ is globally identifiable, i.e., it holds that the matrix
\begin{align*}
H=\frac{1}{n}\sum_{i=1}^n H_i^\top H_i,
\end{align*}
is invertible.
\end{assumption}

The estimation procedure follows an online protocol. At time $t$, each agent $i$ forms an estimate $\te_{i,t}$ of $\theta_t$ based on observations $\{\y_{i,\tau}\}_{\tau=1}^{t-1}$. After that, the new signal $\y_{i,t}$ becomes available to the agent. The online nature of our estimation problem allows us to pose it as an instance of online optimization. Therefore, before deriving the explicit update for $\te_{i,t}$, we need to introduce a few notions based on the literature of online optimization. Let us define a {\it local} square loss 
\begin{align}\label{local loss}
\ell_{i,t}(\theta):=\E\left[\left\|\y_{i,t}-H_i\theta\right\|^2\right],
\end{align}
for each agent $i \in [n]$, followed by the {\it network} loss as 
\begin{align}\label{network loss}
\ell_t(\theta):=\sum_{i=1}^n\E\left[\left\|\y_{i,t}-H_i\theta\right\|^2\right].
\end{align} 
Note that Assumption \ref{A1} guarantees that $\ell_t(\theta)$ has a unique minimizer at each time $t$. 
Agents aim to minimize the accumulated network loss over time. Equivalently, the goal of the network is to minimize the {\it regret} defined as  
\begin{align}\label{regret}
\mathbf{Reg}_T:=\frac{1}{T}\sum_{t=1}^T\bigg(\frac{1}{n} \sum_{j=1}^n \ell_{t}(\te_{j,t})-\ell_t(\theta_t)\bigg).
\end{align}
The regret quantifies the difference between the {\it realized} and {\it ideal} network loss: the estimation $\te_{j,t}$ of each agent $j$ is evaluated at the function $\ell_t(\cdot)$, and the average (among agents) is subtracted by the loss of exact, ideal estimate $\theta_t$. Evidently, the ideal estimate (solution to the network loss \eqref{network loss}) is not available to an individual agent; hence, one can also interpret the regret as the difference between an algorithm that receives local loss (decentralized) versus another algorithm that accesses the network loss (centralized).

\subsection{Network Model and Estimation Update}
The interactions of agents, which in turn defines the network, is captured with the matrix $P$. Formally, we denote by $[P]_{ij}$, the $ij$-th entry of the matrix $P$. When $[P]_{ij}>0$, agent $i$ communicates with agent $j$. We assume that $P$ is symmetric, stochastic with positive diagonal elements. The assumption simply guarantees the information flow in the network. Alternatively, from technical point of view, we respect the following hypothesis. 
\begin{assumption}\label{A2}
We assume the Markov chain $P$ is irreducible and aperiodic, and $P^\top=P$.
\end{assumption}
The assumption implies that the Markov chain $P$ has a unique stationary distribution i.e., $\1^\top P=\1^\top$ is the unique (unnormalized) left eigenvector corresponding to $\lambda_1(P)=1$. It also entails that $\lambda_1(P)$ is unique, and the other eigenvalues of $P$ are less than unit in magnitude \cite{rosenthal1995convergence}. 

To follow the trajectory of the parameter, we propose an update rule which can be cast as a {\it distributed} variant of Online Gradient Descent \cite{zinkevich2003online} with {\it noisy} feedback. It takes the form
\begin{align}\label{update}
\te_{i,t}=\sum_{j=1}^n [P]_{ij} \te_{j,t-1} + \alpha H_i^\top \left(\y_{i,t-1}-H_i\te_{i,t-1}\right),
\end{align}
where $\alpha \in \R$ is the {step size}. The update is akin to {\it consensus+innovation} updates in the literature (see e.g. \cite{khan2010connectivity,kar2012distributed}) with a distinction that it does {\it not} assume any dynamics on the parameter $\theta_t$. The consensus part forces agents to keep their estimates close to each other, and the innovation part takes into account the new observation. We can now shed light on our previous interpretation of regret:  had the agents known the network loss $\ell_t(\theta)$, they could have formed estimates based on Online Gradient Descent (which is a centralized algorithm). Since agents only know the local loss \eqref{local loss}, regret \eqref{regret} is comparing the distributed algorithm versus its centralized counterpart.

We further define the {\it path-length} of the moving target as 
\begin{align}\label{path}
C_T:=\sum_{t=1}^T \| \theta_t-\theta_{t-1}\|^2,
\end{align}
which indicates how much the parameter drifts over time. In particular, when bounding the regret \eqref{regret}, the path-length $C_T$ becomes useful. As we shall explore in the next section, there are connections between mean-square stability and regret; hence, we close this section by defining the mean-square stability. 
\begin{definition}\label{stability}
We say the process is stable in the mean-square sense, whenever
\begin{align*}
\frac{1}{n}\sum_{i=1}^n \E\left[\big\|\te_{i,t}-\theta_t\big\|^2\right] \longrightarrow \sigma,
\end{align*}
where $\sigma\geq0$ is a finite constant.
\end{definition}

\section{Main Results}\label{Main Results}
In this section, we provide our technical results and their consequences. More specifically, our objective is to find an anytime bound on the regret involving the path-length of target. We further interpret the bound in several circumstances. 
\subsection{Error Process and Stability}
Stability (in mean or mean-square sense) is closely attached to the behavior of estimation error. Let us define the {\it local} error process as follows,
\begin{align}\label{local}
\e_{i,t}:=\te_{i,t}-\theta_t.
\end{align}
Stacking the local errors in a vector, we represent the {\it global} error by
\begin{align}\label{global}
\e_t:=[\e_{1,t} , \ldots , \e_{n,t}]^\top.
\end{align}
In the following lemma, we show that the error process can be represented as an LTI system.
\begin{lemma}\label{error lemma}
The error process \eqref{global} can be characterized via an LTI system as follows
\begin{align*}
\e_t=Q \e_{t-1} + \mathbf{u}_t,
\end{align*}
where 
\begin{align*}
Q:=P\otimes I_d-\alpha~\diag\left[H_1^\top H_1,\ldots, H_n^\top H_n\right],
\end{align*}
and 
\begin{align*}
\mathbf{u}_t:= \1_n\otimes(\theta_{t-1}-\theta_t) +\alpha~\left[ \begin{array}{ccc}
H_1^\top \w_{1,t-1} \\
\vdots \\
H_n^\top \w_{n,t-1}  \end{array} \right].
\end{align*}
\end{lemma}
We observe from Lemma \ref{error lemma} that to ensure mean stability (to have $\E[\e_t] \rightarrow \0$ as $t\rightarrow \infty$), the conditions
\begin{align*}
\rho(Q)<1 \ \ \ \ \ \text{and} \ \ \ \ \ \ \E[\theta_t-\theta_{t-1}] \rightarrow \0,
\end{align*}
need to be satisfied. 
To push the spectral radius of $Q$ inside the unit circle, we can simply tune the step size $\alpha$. We quantify this condition in Section \ref{Tuning the Step Size}, and show the role of global identifiability (Assumption \ref{A1}) and connectivity of network (Assumption \ref{A2}) in existence of such $\alpha$.

However, the second condition entirely depends on the problem environment. For instance, the condition is evidently satisfied in the static case ($\theta_t=\theta$). It also holds when $\theta_t$ is generated randomly based on a {\it stationary} distribution, i.e. whenever $\E[\theta_t]=\theta$. Regardless of stability, we can always characterize the regret behavior in terms of path-length of the target.

\subsection{A Path-Length Bound for Regret}\label{A Path-Length Bound for Regret}
In this section, we analyze the regret \eqref{regret}, which is a measure of {\it finite-time} performance. The following lemma exhibits the connection of mean-square stability (Definition \ref{stability}) and regret \eqref{regret}. In particular, the lemma is used to bound the regret in terms of path-length of the target \eqref{path}.

\begin{lemma}\label{regret lemma}
The regret defined in \eqref{regret} can be bounded as,
\end{lemma}
\begin{align*}
\mathbf{Reg}_T \leq \frac{1}{nT} \sum_{t=1}^T\sum_{i=1}^n\left\|H_i\right\|^2\E\left[\left\|\e_t\right\|^2\right].
\end{align*}
Lemma \ref{regret lemma} proves that stability implies boundedness of asymptotic regret since
\begin{align*}
\E\left[\left\|\e_t\right\|^2\right]=\sum_{i=1}^n \E\left[\big\|\te_{i,t}-\theta_t\big\|^2\right],
\end{align*}
and once the limit exists, the Ces{\`a}ro mean preserves it. However, we are interested in non-asymptotic analysis of the problem. We present our finite-time statement in the subsequent theorem. 

\begin{theorem}\label{our theorem}
Let the sequence $\{\theta_t\}_{t=1}^T$ be globally identifiable (Assumption \ref{A1}), and the Markov chain $P$ be irreducible and aperiodic (Assumption \ref{A2}). If each agent $i\in [n]$ generates the estimate sequence $\{\te_{i,t}\}_{t=1}^T$ according to the update rule \eqref{update}, the regret satisfies  
\begin{align*}
\mathbf{Reg}_T &\leq \frac{1}{n} \sum_{i=1}^n\left\|H_i\right\|^2 \frac{\alpha^2 \sum_{i=1}^n\|H_i\|^2W_i}{1-\|Q\|}\\
&~~~~~~~~~~~~~~~~~~+\frac{1}{T} \sum_{i=1}^n\left\|H_i\right\|^2\frac{C_T}{(1-\|Q\|)^2}.
\end{align*}
\end{theorem}
Theorem \ref{our theorem} provides a path-length bound for regret. The underlying intuition behind the term $C_T$ is as follows: as agents decide on the next value based on previous observations (online prediction), they are always one step behind in estimation. Even when the observations are {\it noiseless} ($W_i=0$ for all $i \in [n]$), the second term still remains, which is an artifact of agent $i$ using $\y_{i,t-1}$ to predict $\theta_t$. 

We remark that path-length regret bounds were previously studied in the context of online optimization for {\it centralized} frameworks (see e.g. \cite{zinkevich2003online,hall2015online,jadbabaie2015online}). Here, we specialized to quadratic losses, and proved a path-length bound in {\it distributed} setting. Therefore, our bound involves an additional network penalty comparing to the centralized framework.

\subsection{Tuning the Step Size}\label{Tuning the Step Size}
We now discuss several aspects of the regret bound derived in Theorem \ref{our theorem}. To this end, we need to clarify the dependence of $\|Q\|$ on step size $\alpha$ to analyze the bound. Recall the closed-form of $Q$ in Lemma \ref{error lemma} as well as stochasticity of $P$ which implies $P\1_n=\1_n$; pre/post multiply $Q$ by a vector of all ones with dimension $nd$ ($\1=\1_n\otimes \1_d$) to get
\begin{align*}
\1^\top Q \1 &= \1^\top\left(P\otimes I_d-\alpha~\diag\left[H_1^\top H_1,\ldots, H_n^\top H_n\right]\right)\1\\
&=\1^\top\left(\1-\alpha~\diag\left[H_1^\top H_1,\ldots, H_n^\top H_n\right]\1\right)\\
&=\1^\top\1-\alpha~\1^\top\diag\left[H_1^\top H_1,\ldots, H_n^\top H_n\right]\1\\
&=\1^\top\1-n\alpha~\1_d^\top H\1_d,
\end{align*}
as $H=\frac{1}{n}\sum_{i=1}^nH_i^\top H_i$. Therefore, since $\1^\top\1=nd$, we obtain 
\begin{align*}
\frac{\1^\top Q \1}{\1^\top  \1} &=1-n\alpha~\frac{\1_d^\top H\1_d}{\1^\top\1}=1-\alpha\frac{\1_d^\top H\1_d}{\1_d^\top\1_d}.
\end{align*}
In view of Assumption \ref{A2}, $\lambda_1(P)=1$ is unique, and no direction other than $\1_n$ can recover the trivial eigenvalue. Hence, depending on the null space of $\diag\left[H_1^\top H_1,\ldots, H_n^\top H_n\right]$, $\alpha$ can be set small enough such that the following upper bound holds 
\begin{align}\label{upper}
\lambda_1(Q) \leq  1-\alpha \lambda_n(H).
\end{align}
By the same token, Weyl's eigenvalue inequality implies a lower bound the smallest eigenvalue of $Q$ as
\begin{align}\label{lower}
\lambda_n(Q) \geq \lambda_n(P)-\alpha \lambda_1(\diag\left[H_1^\top H_1,\ldots, H_n^\top H_n\right]). 
\end{align}
Once again choosing small enough $\alpha$ guarantees that the RHS of \eqref{upper} is larger than the absolute value of the RHS of \eqref{lower}. Let us distinguish such regime as $\alpha \leq \alpha_{\max}$. Then, symmetry of $Q$ warrants that
\begin{align*}
\|Q\| \leq 1-\alpha \lambda_n(H), \ \ \ \ \forall \alpha \leq \alpha_{\max}.
\end{align*}
Therefore, disregarding the constants (dependence on $\{H_i\}_{i=1}^n$ and $H$) in Theorem \ref{our theorem}, we can simplify the regret bound using above as follows
\begin{align*}
\mathbf{Reg}_T &\leq \mathcal{O}\left(\alpha \sum_{i=1}^n W_i+ \frac{C_T}{T\alpha^2}\right).
\end{align*}
The following comments are now in order:
\begin{itemize}
\item If the target is fixed, i.e. $\theta_t=\theta$ for every $t\in [T]$, we have $C_T=0$, and the second term in the bound vanishes. Then, we can set $\alpha=\min\{\frac{1}{T},\alpha_{\max}\}$ to maintain the $\mathcal{O}(\frac{1}{T})$ rate. This choice of step size and the corresponding result is consistent with the results in {\it static} setting.
\item When the observations are noiseless ($W_i=0$) the first term becomes zero, and $\alpha=\min\{1,\alpha_{\max}\}$ recovers the $\mathcal{O}(\frac{C_T}{T})$ rate. In this regime, a sub-linear path-length $C_T$ always guarantees a zero asymptotic regret. 
\item In the general case tuning $\alpha=\min\{C_T^{1/3}T^{-1/3},\alpha_{\max}\}$ yields a regret of $\mathcal{O}\left(C_T^{1/3}T^{-1/3}\right)$ which holds for any $T$. 
\end{itemize}

\iffalse
\subsection{Extension to Non-Linear Dynamics}
In the previous section, we considered a moving target that does not follow any particular dynamics. 
We also consider the case that $\theta_t$ follows a non-linear dynamics
\begin{align*}
\theta_{t+1}=f_t(\theta_t),
\end{align*}
where $f_t: \R^n \to \R^n$ is a {\it contractive} mapping for every $t \in [T]$ in the following sense:
\begin{definition}\label{contractive}
A mapping $f: \R^n \to \R^n$ is called contractive when for every $x,y \in \R^n$ it satisfies
\begin{align*}
\|f(x)-f(y)\|\leq \|x-y\|.
\end{align*}
\end{definition}
The knowledge of the dynamics (even if it is non-linear) would help agents to predict the next move of the target.  
\fi

\section{Concluding Remarks}\label{Concluding Remarks}
In this paper, we considered the distributed estimation problem in dynamic environments. A network of agents partially observe the parameter of interest which is time-varying with no particular dynamics. However, the parameter is globally observable from the standpoint of all agents together. We pose the problem as an instance of distributed, online optimization, where agents suffer some loss after each round of estimation. In turn, they aim to minimize the network loss defined as the sum of individual losses. We formulated the problem using the notion of dynamic regret from online optimization literature. Our main contribution was to provide the non-asymptotic analysis of the dynamic regret. We first showed that the regret can be related to the finite-time counterpart of stability. Next, we defined the path-length as the consecutive differences between target locations, and derived an upper bound on regret in terms of the path-length. We finally demonstrated that our bound is consistent with special cases of static setting and noiseless observations.

There are several directions in which our work can be improved. Our setting does not include any prior assumption on the dynamics of the parameter. While the generality of this setup is appealing, it is also interesting to investigate the problem for particular dynamics. More specifically, once agents know the dynamics of the parameter (even in the non-linear case) they can incorporate it into the algorithm, and obtain more accurate estimates. Furthermore, the choice of step size in the general case requires a prior knowledge of the path-length. An alternative approach is to have an online feedback of the trajectory as the algorithm moves forward. Therefore, an adaptive, online algorithm for tracking dynamic parameters is still an open problem.

\section{Appendix : Proofs}\label{Appendix : Proofs}

\noindent
\textbf{\emph{Proof of Lemma \ref{error lemma}}}. In view of equations \eqref{model} and \eqref{update}, we write
\begin{align*}
\te_{i,t}-\theta_t&=\sum_{j=1}^n[P]_{ij}\te_{j,t-1}-\theta_t\\
&+ \alpha H_i^\top \left(\y_{i,t-1}-H_i\te_{i,t-1}\right)\\
&=\sum_{j=1}^n[P]_{ij}(\te_{j,t-1}-\theta_{t-1})+\theta_{t-1}-\theta_t\\
&+ \alpha H_i^\top \left(H_i \theta_{t-1} + \w_{i,t-1}-H_i\te_{i,t-1}\right), \numberthis \label{eq:1}
\end{align*}
since $\sum_{j=1}^n[P]_{ij}=1$. Rewriting above based on the definition of local error process \eqref{local}, we obtain
\begin{align*}
\e_{i,t}&=\sum_{j=1}^n[P]_{ij}\e_{j,t-1}-\alpha H_i^\top H_i\e_{i,t-1} \\
&+\theta_{t-1}-\theta_t+\alpha~H_i^\top \w_{i,t-1}.
\end{align*} 
Representing above in the matrix form using the Kronecker product completes the proof.$\hfill \blacksquare $\\

\noindent
\textbf{\emph{Proof of Lemma \ref{regret lemma}}}. Note equations \eqref{model} and \eqref{network loss} to observe
\begin{align}\label{eq:2}
\ell_t(\theta_t)=\sum_{i=1}^n\E\left[\|\w_{i,t}\|^2\right]=\sum_{i=1}^nW_i,
\end{align}
and 
\begin{align*}
\ell_{t}(\te_{j,t})&=\sum_{i=1}^n\E\left[\left\|\y_{i,t}-H_i\te_{j,t}\right\|^2\right]\\
&=\sum_{i=1}^n\E\left[\left\|H_i\theta_t+\w_{i,t}-H_i\te_{j,t}\right\|^2\right]\\
&=\sum_{i=1}^n\E\left[\left\|H_i\e_{j,t}-\w_{i,t}\right\|^2\right],
\end{align*}
where we used \eqref{local} in the last step. Given that noise realizations are independent over time, we can simplify above
\begin{align*} 
\ell_{t}(\te_{j,t})&=\sum_{i=1}^n\E\left[\left\|H_i\e_{j,t}-\w_{i,t}\right\|^2\right]\\
&=\sum_{i=1}^n\left(\E\left[\left\|H_i\e_{j,t}\right\|^2\right]+W_i\right). \label{eq:3} \numberthis
\end{align*}
Therefore, combining \eqref{eq:2} and \eqref{eq:3}, we get
\begin{align*}
\frac{1}{n} \sum_{j=1}^n \ell_{t}(\te_{j,t})-\ell_t(\theta_t)&=\frac{1}{n} \sum_{j=1}^n\sum_{i=1}^n\E\left[\left\|H_i\e_{j,t}\right\|^2\right]\\
&\leq \frac{1}{n} \sum_{j=1}^n\sum_{i=1}^n\left\|H_i\right\|^2\E\left[\left\|\e_{j,t}\right\|^2\right]\\
&= \Big(\frac{1}{n} \sum_{i=1}^n\left\|H_i\right\|^2\Big)\E\left[\left\|\e_t\right\|^2\right],
\end{align*}
where we used \eqref{global} in the last line. Summing over $t \in [T]$ completes the proof.$\hfill \blacksquare $\\

\noindent
\textbf{\emph{Proof of Theorem \ref{our theorem}}}.  In view of Lemma \ref{error lemma}, we have
\begin{align*}
\E\left[\|\e_t\|^2\right] &\leq \|Q\|^2\E\left[\|\e_{t-1}\|^2\right]\\
&+\alpha^2 \sum_{i=1}^n\|H_i\|^2W_i+n\|\theta_{t-1}-\theta_t\|^2\\
&+2 [\1_n \otimes (\theta_{t-1}-\theta_t)]^\top  Q\E[\e_{t-1}]. \label{eq:4} \numberthis
\end{align*}
since $\E[\w_{i,t-1}]=\0$. Using the simple fact that for any $\beta>0$
\begin{align*}
2ab \leq \frac{1}{\beta}a^2+\beta b^2,
\end{align*}
based on AM-GM inequality, we apply the Cauchy-Schwarz inequality to bound
\begin{align*}
&2 [\1 \otimes (\theta_{t-1}-\theta_t)]^\top  Q\E[\e_{t-1}] \leq \\
&~~~~~~~~~~~~~~~~\frac{1}{\beta} \|\1_n \otimes (\theta_{t-1}-\theta_t)\|^2+\beta\|Q\|^2E[\|\e_{t-1}\|]^2,
\end{align*}
for any $\beta>0$, and simplify \eqref{eq:4} to 
\begin{align*}
\E\left[\|\e_t\|^2\right] &\leq (1+\beta)\|Q\|^2\E\left[\|\e_{t-1}\|^2\right]\\
&+\alpha^2 \sum_{i=1}^n\|H_i\|^2W_i+\frac{\beta+1}{\beta}n\|\theta_{t-1}-\theta_t\|^2,
\end{align*}
which entails
\begin{align*}
&\Big(1-(1+\beta)\|Q\|^2\Big)\E\left[\|\e_t\|^2\right] \\
&~~~~~~~~~~~\leq (1+\beta)\|Q\|^2\Big(\E\left[\|\e_{t-1}\|^2\right]-\E\left[\|\e_t\|^2\right]\Big)\\
&~~~~~~~~~~~+\alpha^2 \sum_{i=1}^n\|H_i\|^2W_i+\frac{\beta+1}{\beta}n\|\theta_{t-1}-\theta_t\|^2. \numberthis \label{eq:5}
\end{align*}
Using the convention $\E\left[\|\e_0\|^2\right]=0$, the following sum telescopes
\begin{align*}
\sum_{t=1}^T\E\left[\|\e_{t-1}\|^2\right]-\E\left[\|\e_t\|^2\right] \leq \E\left[\|\e_0\|^2\right] =0,
\end{align*}
and summing \eqref{eq:5} over $t\in [T]$, we derive
\begin{align*}
\sum_{t=1}^T \E\left[\|\e_t\|^2\right] &\leq \frac{\alpha^2T \sum_{i=1}^n\|H_i\|^2W_i}{1-(1+\beta)\|Q\|^2}\\
&~~~~~~~~~+\frac{\beta+1}{\beta}\frac{n\sum_{t=1}^T\|\theta_{t-1}-\theta_t\|^2}{1-(1+\beta)\|Q\|^2}.
\end{align*}
Recall that $0<\beta<\|Q\|^{-2}-1$ is a free parameter in the analysis. Letting $\beta=\|Q\|^{-1}-1$, and recalling Definition \ref{path}, we calculate the bound above to get
\begin{align*}
\sum_{t=1}^T \E\left[\|\e_t\|^2\right] &\leq \frac{\alpha^2T \sum_{i=1}^n\|H_i\|^2W_i}{1-\|Q\|}+\frac{nC_T}{(1-\|Q\|)^2}.
\end{align*}
Appealing to Lemma \ref{regret lemma} completes the proof.$\hfill \blacksquare $\\

\bibliographystyle{IEEEtran}
\bibliography{IEEEabrv,shahin}

\end{document}